\documentclass[a4paper,12pt]{article}
\usepackage[ansinew]{inputenc}
\usepackage{amsfonts}
\usepackage{latexsym}
\usepackage{amsmath}
\usepackage{amssymb}
\usepackage{jeep}

\newcommand{\R}{\mathbb{R}}
\newcommand{\C}{\mathbb{C}}
\newcommand{\N}{\mathbb{N}}
\newcommand{\Z}{\mathbb{Z}}

\newtheorem{defin}{Definition}[section]
\newtheorem{theorem}[defin]{Theorem}
\newtheorem{lemma}[defin]{Lemma}
\newtheorem{corollary}[defin]{Corollary}
\newenvironment{proof}
{\noindent{\it Proof.}}{\hfill $\Box$\par\vspace{2.5mm}}

\numberwithin{equation}{section}

\title{Nevanlinna theory for the difference operator}
\author{R. G. Halburd and R. J. Korhonen}
\date{}

\begin{document}

\maketitle

\begin{abstract}
Certain estimates involving the derivative $f\mapsto f'$ of a
meromorphic function play key roles in the construction and
applications of classical Nevanlinna theory. The purpose of this
study is to extend the usual Nevanlinna theory to a theory for the
exact difference $f\mapsto \Delta f=f(z+c)-f(z)$.

An $a$-point of a meromorphic function $f$ is said to be
$c$-\textit{paired} at $z\in\C$ if $f(z)=a=f(z+c)$ for a fixed
constant $c\in\C$. In this paper the distribution of paired points
of finite-order meromorphic functions is studied. One of the main
results is an analogue of the second main theorem of Nevanlinna
theory, where the usual ramification term is replaced by a
quantity expressed in terms of the number of paired points of $f$.
Corollaries of the theorem include analogues of the Nevanlinna
defect relation, Picard's theorem and Nevanlinna's five value
theorem. Applications to difference equations are discussed, and a
number of examples illustrating the use and sharpness of the
results are given.
\end{abstract}

\section{Introduction}

\renewcommand{\thefootnote}{}
\footnotetext[1]{\emph{Mathematics Subject Classification 2000:}
Primary 30D35; Secondary 39A10, 39A12.} \footnotetext[2]{The
research reported in this paper was supported in part by the
Leverhulme Trust and by the Finnish Academy (grant number
204819).}

Nevanlinna's theory of value distribution is concerned with the
density of points where a meromorphic function takes a certain
value in the complex plane. One of the early results in this area
is a theorem by Picard \cite{picard:80} which states that a
non-constant entire function can omit at most one value.
Nevanlinna offered a deep generalization of Picard's theorem in
the form of his \textit{second main theorem} \cite{nevanlinna:25},
which implies the defect relation:
    \begin{equation}\label{defsum}
    \sum_a \left(\delta(a,f) +  \theta(a,f)\right) \leq 2
    \end{equation}
where the sum is taken over all points in the extended complex
plane, $f$ is a non-constant meromorphic function and the
quantities $\delta(a,f)$ and $\theta(a,f)$ are called the
\textit{deficiency} and the \textit{index of multiplicity} of the
value $a$, respectively (see Section~\ref{notation}). The defect
relation \eqref{defsum} yields, for instance, Picard's theorem as
an immediate corollary. It also implies that the maximum number of
totally ramified values is at most four for any meromorphic
function.

The appearance of the ramification index $\theta(a,f)$ in the
defect relation \eqref{defsum} means that the density of
high-multiplicity $a$-points is relatively low for most $a\in\C$.
Similarly in this paper it is shown that $a$-points appearing in
pairs with constant separation are rare for finite-order
meromorphic functions, unless the function in question is periodic
with the same period as the separation. For instance, if $f$ is of
finite order and not periodic with period $c$, then
    \begin{equation}\label{defsumpair}
    \sum_a \left(\delta(a,f) +  \pi_c(a,f)\right) \leq 2
    \end{equation}
where the sum is taken over all points in the extended complex
plane, and $\pi_c(a,f)$ is a measure of those $a$-points of $f$
which appear in pairs separated by the constant $c\in\C$ (in other
words, those points $z_0$ where $f(z_0)=a=f(z_0+c)$, see
Section~\ref{notation} for the exact definition.) The sharpness of
inequality \eqref{defsumpair} is shown by giving an example of a
finite-order meromorphic function, which is not periodic with
period~$c$, satisfying $\sum_a\pi_c(a,f)=2$.

The defect relation \eqref{defsum} follows by an analysis of the
behavior of the derivative $f\mapsto f'$ in the ramification term
of the second main theorem. In what follows, \eqref{defsumpair} is
obtained by proving a version of the second main theorem where the
derivative of $f$ is replaced by the exact difference
$f\mapsto\Delta f=f(z+c)-f(z)$ of a meromorphic function. In the
remainder of this paper difference analogues of Picard's theorem
and Nevanlinna's theorem on functions sharing five values are
given. In addition, the sharpness of the obtained results is
discussed with the help of examples, and an application to
difference equations is presented.

\section{Nevanlinna theory for exact differences}

Before going into details of value distribution of exact
differences we must first give a precise answer to the following
question: What is the difference analogue of a point with high
multiplicity? By a formal discretisation of the derivative
function $f'(z)$ we obtain
    \begin{equation}\label{discr}
    \frac{f(z+c)-f(z)}{c}=:\frac{\Delta_c f}{c},
    \end{equation}
where $c\in\C$. As noted in the introduction, those $a$-points of
$f$ where the derivative vanishes, called \textit{ramified}
points, play a special role in Nevanlinna theory. The
discretisation \eqref{discr} of $f'(z)$ suggests that $a$-points
appearing in pairs separated by a fixed constant $c$ may have
similar importance with respect to the operator $\Delta_c$. This
indeed turns out to be the case as seen in the following sections.

\subsection{Lemma on the exact difference}\label{notation}

We first briefly recall some of the basic definitions of
Nevanlinna theory. We refer to \cite{hayman:64} for a
comprehensive description of the value distribution theory. The
Nevanlinna \textit{deficiency} is defined as
    $$
    \delta(a,f) := \liminf_{r\to\infty} \frac{m(r,a)}{T(r,f)},
    $$
where $a$ is in the extended complex plane, $m(r,a)$ is the
Nevanlinna \textit{proximity function} and $T(r,f)$ is the
\textit{characteristic function} of $f$. The \textit{ramification
index} is
    $$
    \theta(a,f) := \liminf_{r\to\infty}
    \frac{N(r,a)-\overline{N}(r,a)}{T(r,f)},
    $$
where $N(r,a)$ is the \textit{counting function} of the $a$-points
of $f$, counting multiplicities, and $\overline{N}(r,a)$ the
counting function ignoring multiplicities. The point $a\in\C$ is a
\textit{totally ramified value} of $f$ if all $a$-points of $f$
have multiplicity two or higher.

The following theorem is a recently obtained difference analogue
of the lemma on the logarithmic derivative \cite{halburdk:prep04}.

\begin{theorem}\label{logdiff}
Let $f$ be a non-constant meromorphic function of finite order,
$c\in\C$ and $\delta<1$. Then
    \begin{equation}\label{diff}
    m\left(r,\frac{f(z+c)}{f(z)}\right)= o\left(\frac{T(r,f)}{r^\delta}\right)
    \end{equation}
for all $r$ outside of a possible exceptional set $E$ with finite
logarithmic measure $\int_E\frac{dr}{r}<\infty$.
\end{theorem}

In the original statement of Theorem~\ref{logdiff} in
\cite{halburdk:prep04} the error term on the right side of
\eqref{diff} has $T(r+|c|,f)$ instead of $T(r,f)$. But by the
following lemma, \cite[Lemma~2.1]{halburdk:04}, we have
$T(r+|c|,f)=(1+o(1))T(r,f)$ for all $r$ outside of a set with
finite logarithmic measure, whenever $f$ is of finite order.

\begin{lemma}\label{technical}
Let $T:(0,+\infty)\to(0,+\infty)$ be a non-decreasing continuous
function, $s>0$, $\alpha<1$, and let $F\subset\R^{+}$ be the set
of all $r$ such that
    \begin{equation}\label{assu}
    T(r) \leq \alpha T(r+s).
    \end{equation}
If the logarithmic measure of is $F$ infinite, that is,
$\int_F\frac{dt}{t}=\infty$, then
    $$
    \limsup_{r\to\infty}\frac{\log T(r)}{\log r}=\infty.
    $$
\end{lemma}

Let $f(z)$ be a non-constant meromorphic function of finite order,
and let $a(z)$ be a finite-order periodic function with period $c$
such that $f(z)\not\equiv a(z)$. Denote
    $$
    \Delta_c f := f(z+c) - f(z),
    $$
and $\Delta^n_c f:=\Delta^{n-1}_c( \Delta_c f)$ for all $n\in\N$,
$n\geq2$. Then by applying Theorem~\ref{logdiff} with the function
$f(z)-a(z)$, we have
    \begin{equation}\label{logdiffeq2}
    \begin{split}
    m\left(r,\frac{\Delta_c f}{f-a}\right) &= m\left(r,\frac{f(z+c)-a(z+c)}{f(z)-a(z)}\right)  + O(1)\\
    & =o\left(\frac{T(r,f-a)}{r^\delta}\right) + O(1)
    \end{split}
    \end{equation}
outside of a possible exceptional set with finite logarithmic
measure. We denote by $\mathcal{S}(f)$ the set of all  meromorphic
functions $g$ such that $T(r,g)=o(T(r,f))$ for all $r$ outside of
a set with finite logarithmic measure. Functions in the set
$\mathcal{S}(f)$ are called \textit{small compared to} $f$, or
\textit{slowly moving} with respect to $f$. Also, if
$g\in\mathcal{S}(w)$ we denote $T(r,g)=S(r,f)$.

Since by \eqref{logdiffeq2}
    \begin{equation}\label{logdiffeq}
    m\left(r,\frac{\Delta_c f}{f-a}\right) =S(r,f-a)
    \end{equation}
we arrive at the following lemma by induction and using the fact
that
    \begin{equation*}
    T(r,f(z+1))\leq (1+\varepsilon)T(r+1,f(z))
    \end{equation*}
for any $\varepsilon>0$ when $r$ is large \cite{yanagihara:80}.

\begin{lemma}
Let $c\in\C$, $n\in\N$, and let $f$ be a meromorphic function of
finite order. Then for all small periodic functions
$a\in\mathcal{S}(f)$
    \begin{equation*}
    m\left(r,\frac{\Delta^n_c f}{f-a}\right) =S(r,f),
    \end{equation*}
where the exceptional set associated with $S(r,f)$ is of at most
finite logarithmic measure.
\end{lemma}

Finally, an identity due to Valiron \cite{valiron:31} and Mohon'ko
\cite{mohonko:71} is needed in the following section. It states
that if the function $R(z,f)$ is rational in $f$ and has small
meromorphic coefficients, then
    \begin{equation}\label{vm}
    T(r,R(z,f)) = \deg_f(R)T(r,f) + S(r,f).
    \end{equation}
For the proof see also~\cite{laine:93}.

\subsection{Second main theorem}

The lemma on the logarithmic derivative is one of the main
components of the proof of the second main theorem of Nevanlinna
theory. The following theorem is obtained by combining the
standard method of proof for the second main theorem
\cite{nevanlinna:25} together with Theorem~\ref{logdiff}. As a
result a version of the second main theorem is obtained where,
instead of the usual ramification term, there is a certain
quantity expressed in terms of paired points of the considered
function $f$. Since periodic functions are the analogues of
constants for exact differences, it is natural to consider slowly
moving periodic functions as target functions of $f$.

\begin{theorem}\label{2nd}
Let $c\in\C$, and let $f$ be a meromorphic function of finite
order such that $\Delta_c f\not\equiv 0$. Let $q\geq 2$, and let
$a_1(z),\ldots,a_q(z)$ be distinct meromorphic periodic functions
with period $c$ such that $a_k\in \mathcal{S}(f)$ for all
$k=1,\ldots,q$. Then
    \begin{equation*}
    m(r,f) + \sum_{k=1}^q m\left(r,\frac{1}{f-a_k}\right) \leq 2T(r,f) -N_{pair}(r,f)
     + S(r,f)
    \end{equation*}
where
    $$
    N_{pair}(r,f):=2N(r,f)-N(r,\Delta_c f)+N\left(r,\frac{1}{\Delta_c f}\right)
    $$
and the exceptional set associated with $S(r,f)$ is of at most finite logarithmic measure.
\end{theorem}

\begin{proof}
By denoting
    $$
    P(f):= \prod_{k=1}^q \left(f-a_k\right),
    $$
we have
    $$
    \frac{1}{P(f)}=\sum_{k=1}^q\frac{\alpha_k}{f-a_k},
    $$
where $\alpha_k\in\mathcal{S}(f)$ are certain periodic functions
with period $c$. Hence, by \eqref{logdiffeq}, we obtain
    $$
    m\left(r,\frac{\Delta_c f}{P(f)}\right) \leq \sum_{k=1}^q m\left(r,\frac{\Delta_c f}{f-a_k}\right) +
    S(r,f)=S(r,f),
    $$
and so
    \begin{equation}\label{estim}
    m\left(r,\frac{1}{P(f)}\right) = m\left(r,\frac{\Delta_c f}{P(f)}\frac{1}{\Delta_c f}\right) \leq
    m\left(r,\frac{1}{\Delta_c f}\right) + S(r,f).
    \end{equation}
By combining the first main theorem, \eqref{estim} and the Valiron-Mo'honko identity \eqref{vm}, we have
    \begin{equation*}
    \begin{split}
    T(r, \Delta_c f) & = m\left(r,\frac{1}{\Delta_c f}\right) + N\left(r,\frac{1}{\Delta_c f}\right) + O(1)\\
    & \geq m\left(r,\frac{1}{P(f)}\right) + N\left(r,\frac{1}{\Delta_c f}\right) +S(r,f) \\
    & = qT(r,f) - \sum_{k=1}^q N\left(r,\frac{1}{f-a_k}\right) + N\left(r,\frac{1}{\Delta_c f}\right)
    +S(r,f)\\
    &= \sum_{k=1}^q m\left(r,\frac{1}{f-a_k}\right) + N\left(r,\frac{1}{\Delta_c f}\right)
    +S(r,f).\\
    \end{split}
    \end{equation*}
Thus, by \eqref{logdiffeq},
    \begin{equation*}
    \begin{split}
    m(r,f) + \sum_{k=1}^q m\left(r,\frac{1}{f-a_k}\right)
    &\leq T(r,f) +N(r,\Delta_c f) + m(r,\Delta_c f) \\ &\quad -
    N\left(r,\frac{1}{\Delta_c f}\right) - N(r,f) +S(r,f)\\
    & \leq T(r,f) +N(r,\Delta_c f) + m(r,f) \\ &\quad -
    N\left(r,\frac{1}{\Delta_c f}\right) - N(r,f) +S(r,f)\\
    &=2T(r,f) +N(r,\Delta_c f)   -
    N\left(r,\frac{1}{\Delta_c f}\right)\\ &\quad - 2N(r,f) +S(r,f).
    \end{split}
    \end{equation*}
\end{proof}

Let us now analyze the assertion of Theorem~\ref{2nd} more
closely. By Lemma \ref{technical} $N(r+|c|,f)=(1+o(1))N(r,f)$ for
all $r$ outside of a set with finite logarithmic measure.
Therefore,
    \begin{equation*}
    \begin{split}
    N_{pair}(r,f) &\geq N(r,f) - N(r+|c|,f) + N\left(r,\frac{1}{\Delta_c f}\right)\\
    &=N\left(r,\frac{1}{\Delta_c f}\right)+
     S(r,f)
    \end{split}
    \end{equation*}
so clearly Theorem~\ref{2nd} is telling us something non-trivial
about the value distribution of finite-order meromorphic
functions. In order to better interpret the meaning of the
\textit{pair term} $N_{pair}(r,f)$ we introduce the counting
function $n_c(r,a)$, $a\in\C$, which is the number of points $z_0$
where $f(z_0)=a$ and $f(z_0+c)=a$, counted according to the number
of equal terms in the beginning of Taylor series expansions of
$f(z)$ and $f(z+c)$ in a neighborhood of $z_0$. We call such
points $c$-\textit{separated} $a$-\textit{pairs} of $f$ in the
disc $\{z:|z|\leq r\}$.

For instance, if $f(z)=a$ and $f(z+c)=a$ with multiplicities $p$
and $q<p$, respectively, then the $q$ first terms in the series
expansions of $f(z)$ and $f(z+c)$ are identical, and so this point
is counted $q$ times in $n_c(r,a)$. Similarly, if in a
neighborhood of $z_0$
    $$
    f(z)=a+c_1(z-z_0)+c_2(z-z_0)^2 + \alpha(z-z_0)^3
    +O\left((z-z_0)^4\right)
    $$
and
    $$
    f(z+c)=a+c_1(z-z_0)+c_2(z-z_0)^2 + \beta(z-z_0)^3
    +O\left((z-z_0)^4\right)
    $$
where $\alpha\not=\beta$, then the point $z_0$ is counted $3$
times in $n_c(r,a)$.

The integrated counting function is defined as follows:
    $$
    N_c(r,a):=\int_0^r\frac{n_c(t,a)-n_c(0,a)}{t}\,dt + n_c(0,a)\log
    r.
    $$
Similarly,
    $$
    N_c(r,\infty):=\int_0^r\frac{n_c(t,\infty)-n_c(0,\infty)}{t}\,dt
     + n_c(0,\infty)\log r,
    $$
where $n_c(r,\infty)$ is the number of $c$-\textit{separated pole
pairs} of $f$, which are exactly the $c$-separated $0$-pairs of
$1/f$. This means that if $f$ has a pole with multiplicity $p$ at
$z_0$ and another pole with multiplicity $q$ at $z_0+c$ then this
pair is counted $\min\{p,q\}+m$ times in $n_c(r,\infty)$, where
$m$ is the number of equal terms in the beginning of the Laurent
series expansions of $f(z)$ and $f(z+c)$ in a neighborhood of
$z_0$. Of course, if $p\not=q$ then $m=0$.

Note that $n_c(r,a)$ is finite for any finite $r$, provided that
the given function $f$ is not periodic with period $c$. Otherwise
there would be a point $z_0\in\C$ in a neighborhood of which the
series expansions of $f(z)$ and $f(z+c)$ would be identical. But
this means that $f(z)\equiv f(z+c)$ in the whole complex plane,
which contradicts the assumption. However, it is possible that
$n_c(r,a)$ is strictly greater than the counting function
$n(r,a)$.

A natural difference analogue of $\overline{N}(r,a)$ is
    $$
    \widetilde{N}_c(r,a):= N(r,a)-N_c(r,a)
    $$
which counts the number of those $a$-points (or poles) of $f$
which are \textit{not} in $c$-separated pairs. We also use the
notation $N_c(r,\frac{1}{f-a})$ instead of $N_c(r,a)$ and
$N_c(r,f)$ instead of $N_c(r,\infty)$ when we want to emphasize
the connection to the meromorphic function $f$. With this notation
we may state the main result of this paper.

\begin{theorem}\label{2nd2}
Let $c\in\C$, and let $f$ be a meromorphic function of finite
order such that $\Delta_c f\not\equiv 0$. Let $q\geq 2$, and let
$a_1(z),\ldots,a_q(z)$ be distinct meromorphic periodic functions
with period $c$ such that $a_k\in \mathcal{S}(f)$ for all
$k=1,\ldots,q$. Then
    \begin{equation*}\label{TNineq}
    (q-1)T(r,f)  \leq \widetilde{N}_c(r,f) + \sum_{k=1}^q \widetilde{N}_c\left(r,\frac{1}{f-a_k}\right) + S(r,f)
    \end{equation*}
where the exceptional set associated with $S(r,f)$ is of at most
finite logarithmic measure.
\end{theorem}

Before proving Theorem \ref{2nd2} we briefly discuss its
implications. Analogously to the classical Nevanlinna theory, the
counting function $\widetilde{N}_c(r,a)$ satisfies
$\widetilde{N}_c(r,a)=T(r,f)+S(r,f)$ for all except at most
countably many values $a$ (see \cite[pp. 43-44]{hayman:64} for a
proof of this). However, unlike $N(r,a)$, the counting function
$\widetilde{N}_c(r,a)$ may, for some values $a$, be negative for
all sufficiently large $r$. This fact has interesting
consequences. By Theorem~\ref{2nd2} any finite-order meromorphic
function $f$ is either periodic with period $c$, or it can have at
most one non-deficient value $a$ such that whenever $f(z)=a$ also
$f(z+c)=a$ and the first two terms in the series expansions of
$f(z)$ at $z$ and $z+c$ are identical. For instance, consider the
function $g(z):=\wp(z)+\exp(z)$ where $\wp(z)$ is a Weierstrass
elliptic function with a period $c\not=2\pi i$. Then
$T(r,g)=N(r,g)+S(r,g)$ and each pole of $g$ contributes $2$ to
$n(r,g)$ but $-2$ to $\widetilde{n}_c(r,g)$. Therefore
$T(r,g)=-\widetilde{N}_c(r,g)+S(r,g)$ and so
$\widetilde{N}_c(r,a)=T(r,g)+S(r,g)$ for all $a\in\C$ by Theorem
\ref{2nd2}.

\bigskip

\noindent{\it Proof of Theorem \ref{2nd2}.} By Theorem~\ref{2nd}
and the first main theorem, we obtain
    \begin{equation}\label{start}
    \begin{split}
    (q-1)T(r,f)  &\leq N(r,f) +  \sum_{k=1}^q N\left(r,\frac{1}{f-a_k}\right)
    -N\left(r,\frac{1}{\Delta_c f}\right)\\
     &\quad + N(r,\Delta_c f)-2N(r,f)+ S(r,f).
    \end{split}
    \end{equation}
We denote by $N_0(r,f)$ the counting function for those poles of
$f$ having Laurent series expansions at $z_0$ and $z_0+c$ with
identical principal parts, multiplicity counted according to the
number of equal terms in the beginning of the analytic part of the
series expansions. (For instance, if $f(z)=c/(z-z_0)^2
+b/(z-z_0)+a+\alpha(z-z_0)+O((z-z_0)^2)$ and $f(z+c)=c/(z-z_0)^2
+b/(z-z_0)+a+\beta(z-z_0)+O((z-z_0)^2)$ the pole at $z_0$ is
counted once in $N_0(r,f)$ whenever $\alpha\not=\beta$.) Since
$N(r,f)=N(r+|c|,f)+S(r,f)$ by Lemma~\ref{technical}, inequality
\eqref{start} takes the form
    \begin{equation}\label{start2}
    \begin{split}
    (q-1)T(r,f)  &\leq N(r,f) + N_0(r,f)+  \sum_{k=1}^q N\left(r,\frac{1}{f-a_k}\right)
    -N\left(r,\frac{1}{\Delta_c f}\right)\\
     &\quad + N(r,\Delta_c f)-2N(r+|c|,f)-N_0(r,f)+ S(r,f).
    \end{split}
    \end{equation}

The rest of the proof consists of estimates on different terms on
the right side of \eqref{start2}. First, by the definition of a
paired point, we have
    \begin{equation*}
    N_0(r,f)+\sum_{k=1}^q N_c\left(r,\frac{1}{f-a_k}\right) \leq N\left(r,\frac{1}{\Delta_c f}\right)
    \end{equation*}
for all $r$, and thus
    \begin{equation}\label{others}
    N_0(r,f)+\sum_{k=1}^q N\left(r,\frac{1}{f-a_k}\right)-N\left(r,\frac{1}{\Delta_c f}\right)
     \leq \sum_{k=1}^q
     \widetilde{N}_c\left(r,\frac{1}{f-a_k}\right).
    \end{equation}
Second, assume that $z_0\in\C$ is such that $f(z_0+kc)=\infty$ for
all $k\in\Z$ with multiplicities $p_k\geq 0$. Here $p_{k}=0$ means
that $f(z_0+kc)$ is finite. (Note that the case $p_k=0$ for all
$k\not=0$ is not ruled out.) Out of these points only finitely
many are inside the disc $\{z\in\C:|z|\leq r+|c|\}$ for any $r>0$.
By redefining $z_0$ if necessary, we may assume that these points
are $z_0+jc$, $j=0,\ldots,K$, where $K\in\N$ is a constant
depending only on $r$. Then $z_0+c,\ldots,z_0+(K-1)c$ are inside
the disc with radius $r$ centered at the origin, and $\Delta_c f$
has a pole with multiplicity $\max\{p_j,p_{j+1}\}-m_j'$ at
$z_0+jc$, where $j=1,\ldots,K-1$ and $m_j'$ is the number of equal
terms in the beginning of the principal parts of the Laurent
series expansions of $f(z)$ and $f(z+c)$ at $z_0+jc$. If principal
parts are completely identical, the number of equal terms in the
beginning of the analytic parts of the series at $z_0+jc$ is
denoted by $m_j''$, and moreover $m_j:=m_j'+m_j''$. Therefore the
contribution to
    $$
    n(r,\Delta_c f)-2n(r+|c|,f)-n_0(r,f)
    $$
from the points $z_0+jc$, $j=0,\ldots,K$, is
    \begin{equation}\label{2poles}
    \begin{split}
    \sum_{j=1}^{K-1} &  \left( \max\{p_j,p_{j+1}\}-m_j'\right)-2\sum_{j=0}^K
    p_j -\sum_{j=1}^{K-1} m_j''\\
    &= \sum_{j=1}^{K-1} \left( \max\{p_j,p_{j+1}\}-m_j'-m_j''\right) \\ &\quad -\left(p_0 + \sum_{j=0}^{K-1}
    \left(\max\{p_j,p_{j+1}\}+\min\{p_j,p_{j+1}\}\right)+p_K\right)\\
    &\le -\sum_{j=1}^{K-1} \left(\min\{p_j,p_{j+1}\}+m_j\right).
    \end{split}
    \end{equation}
The quantity on the right side of \eqref{2poles} is by definition
exactly the same as the contribution to $-n_c(r,f)$ from the
points $z_0+jc$, $j=0,\ldots,K$. Therefore, by summing over all
poles of $f$, we obtain
    \begin{equation}\label{ineq}
    N(r,f)+N(r,\Delta_c f)-2N(r+|c|,f)-N_0(r,f) \leq
    \widetilde{N}_c(r,f).
    \end{equation}
The assertion follows by combining \eqref{start2}, \eqref{others}
and \eqref{ineq}. \hfill $\Box$\par\vspace{2.5mm}

\subsection{Defect relation and Picard's theorem}

Nevanlinna's second main theorem is a deep generalization of
Picard's theorem, and as such it has many important consequences
for the value distribution of meromorphic functions. In this
section we present difference analogues of a number of these
results, including Picard's theorem and Nevanlinna's theorems on
the total deficiency sum and completely ramified values of a
meromorphic function. All of the results in this section follow
from Theorem~\ref{2nd2}.

A difference analogue of the index of multiplicity $\theta(a,f)$
is called the $c$-\textit{separated pair index}, and it is defined
as follows:
    $$
    \pi_c(a,f):= \liminf_{r\to \infty} \frac{N_c(r,a)}{T(r,f)},
    $$
where $a$ is either a slowly moving periodic function with period
$c$, or $a=\infty$. Similarly, we define
    $$
    \Pi_c(a,f):= 1-\limsup_{r\to \infty}
    \frac{\widetilde{N}_c(r,a)}{T(r,f)},
    $$
which is an analogue of
    $$
    \Theta(a,f)=1-\limsup_{r\to\infty} \frac{\overline{N}(r,a)}{T(r,f)}
    $$
in the usual value distribution theory.

The following corollary says that a non-periodic meromorphic
function of finite order cannot have too many $a$-points which
appear in pairs. It is a difference analogue of Nevanlinna's
theorem on deficient values.

\begin{corollary}\label{rel}
Let $c\in\C$, and let $f$ be a meromorphic function of finite
order such that $\Delta_c f\not\equiv 0$. Then $\Pi_c(a,f)=0$
except for at most countably many meromorphic periodic functions
$a$ with period $c$ such that $a\in \mathcal{S}(f)$, and
    \begin{equation}\label{relsum}
    \sum_{a} \left(\delta(a,f)+\pi_c(a,f)\right)\leq \sum_{a} \Pi_c(a,f) \leq 2.
    \end{equation}
\end{corollary}

By the second main theorem it follows that $\Theta(a,f)=0$ for all
except at most countably many values $a$, see, for instance,
\cite[pp. 43--44]{hayman:64}. The same reasoning can be applied to
prove that Theorem \ref{2nd2} implies Corollary \ref{rel}.

Probably the most distinct difference between the classical
Nevanlinna theory and its difference analogue is that, although
$0\leq \Theta(a,f)\leq 1$ for all meromorphic functions $f$ and
for all $a$ in the extended complex plane, the maximal deficiency
sum
    $$
    \sum_{a} \Pi_c(a,f) = 2
    $$
may be attained by a single value $a$. For instance, the function
$g(z)=\wp(z)+\exp(z)$, where $\wp(z)$ is a Weierstrass elliptic
function with a period $c\not=2\pi i$, satisfies
$\Pi_c(\infty,g)=2$. In fact, by the definition of $\Pi_c(a,f)$
alone, it is not even clear that $\Pi_c(a,f)$ has an upper bound
what so ever. The fact that $\Pi_c(a,f)\leq2$ for all $a$ follows
by Corollary~\ref{rel}.

We say that $a$ is an \textit{exceptional paired value} of $f$
\textit{with the separation} $c$ if the following property holds
for all except at most finitely many $a$-points of $f$:
\textit{Whenever $f(z)=a$ then also $f(z+c)=a$ with the same or
higher multiplicity.} Clearly $N(r,a)\le N_c(r,a)$ for all
exceptional paired values $a$ of $f$. Note also that by this
definition all Picard exceptional values of $f$ are also
exceptional paired values. The following corollary is an analogue
of Picard's theorem.

\begin{corollary}\label{picard}
If a finite-order meromorphic function $f$ has three exceptional
paired values with the separation $c$, then $f$ is a periodic
function with period~$c$.
\end{corollary}

Corollary \ref{picard} implies that if a finite-order meromorphic
function $w$ has two groups of three exceptional paired values
with two different separations, say $c_1$ and $c_2$ independent
over the reals, then either $w$ is a constant or $w$ is an
elliptic function with periods $c_1$ and $c_2$ and therefore
exactly of order $2$.

There is no hope of extending Corollary~\ref{picard} (or Corollary
\ref{rel}) to include all infinite order meromorphic functions,
since the function $\exp(\exp(z))$ has three exceptional paired
values with the separation $\log 2$: In addition to the Picard
exceptional zeros and poles, the value $1$ is exceptionally
paired, although non-deficient.

An example of a finite-order meromorphic function which has
exactly two exceptional paired values with the separation $2K$ is
given by the elliptic function $\textrm{sn}(z,k)$, where
$k\in(0,1)$ is the elliptic modulus and $K$ is the complete
elliptic integral. The function $\textrm{sn}(z,k)$ is periodic
with the periods $4K$ and $2iK'$, and it attains the value zero at
points $2nK+2miK'$ and has its poles at $2nK+(2m+1)iK'$, where
$n,m\in\Z$. The function $\textrm{sn}(z,k)$ has no deficient
values, but it has the maximal four completely ramified values at
$\pm1$ and $\pm1/k$. Therefore, the function
$g(z)=\textrm{sn}(z,k)$ satisfies
    $$
    \sum_a\pi_{2K}(a,g)=2
    $$
and, moreover,
    $$
    \sum_a \left(\theta(a,g)+ \pi_{2K}(a,g)\right)=4.
    $$

Analogously to complete ramification, we say that a point $a$ is
\textit{completely paired with the separation} $c$  if whenever
$f(z)=a$ then either $f(z+c)=a_j$ or $f(z-c)=a_j$, with the same
multiplicity. Then a non-periodic meromorphic function of finite
order can have at most four values which only appear in pairs.

\begin{corollary}
Let $c\in\C$, and let $f$ be a meromorphic function of finite
order such that $\Delta_c f\not\equiv 0$. Then $f$ has at most
four completely paired points with separation $c$.
\end{corollary}

Similarly, a non-periodic finite-order function $f$ can have at
most three values~$a$ which only appear such that for some
$z_0\in\C$, $f(z_0)\not=a$, $f(z_0+jc)=a$ with the same
multiplicity for each $j=1,2,3$, and $f(z_0+4c)\not=a$. We say
that such values appear in lines of three. Similarly, a
finite-order meromorphic function can have a maximum of two values
which appear only in lines of four or more.

\subsection{Functions sharing values}

Another consequence of Nevanlinna's second main theorem is the
five value theorem, which says that if two non-constant
meromorphic functions share five values ignoring multiplicity then
these functions must be identical. By considering periodic
functions instead of constants, and by ignoring paired points
instead of multiplicity, we obtain a difference analogue of the
five value theorem.

We say that two meromorphic functions $f$ and $g$ share a point
$a$, \textit{ignoring} $c$-\textit{separated pairs}, when $f(z)=a$
if and only if $g(z)=a$ with the same multiplicity, unless $a$ is
a $c$-separated pair of $f$ or $g$. In short, all paired points
are ignored when determining whether or not $f$ and $g$ share $a$.
This also means that if $f$ has a paired $a$-point at $z_0$ and
$g$ has a single $a$-point at the same location, this point is not
shared by $f$ and $g$.

\begin{theorem}\label{5}
Let $c\in\C$, and let $f$ and $g$ be meromorphic functions of
finite order. If there are five distinct periodic functions
$a_k\in\mathcal{S}(f)$ such that $f$ and $g$ share $a_k$, ignoring
$c$-separated pairs, for all $k=1,\ldots,5$ then either
$f(z)\equiv g(z)$ or both $f$ and $g$ are periodic with period
$c$.
\end{theorem}

\begin{proof}
We follow the reasoning of the proof of the five value theorem
\cite{hayman:64}. Suppose first that $f$ is periodic with period
$c$. Then by definition all $a$-points of $f$ are paired. Since
$f$ and $g$ share five points, ignoring pairs, $g$ has at least
five exceptional paired values, and therefore it must also be
periodic with period $c$ by Corollary~\ref{picard}.

Assume now that neither $f$ nor $g$ is periodic with period $c$,
and that  $f\not\equiv g$. Then by Theorem~\ref{2nd2}, for any
$\varepsilon>0$,
    \begin{equation}\label{f}
    (4+\varepsilon)T(r,f) \leq \widetilde{N}_c(r,f) +  \sum_{k=1}^5 \widetilde{N}_c\left(r,\frac{1}{f-
    a_k}\right)
    \end{equation}
and
    \begin{equation}\label{g}
    (4+\varepsilon)T(r,g) \leq \widetilde{N}_c(r,g) +  \sum_{k=1}^5 \widetilde{N}_c\left(r,\frac{1}{g-
    a_k}\right)
    \end{equation}
outside a set with finite logarithmic measure.  Since
    $$
    \widetilde{N}_c\left(r,\frac{1}{g-  a_k}\right) = \widetilde{N}_c\left(r,\frac{1}{f-
    a_k}\right)
    $$
for all $k=1,\ldots,5$, inequalities \eqref{f} and \eqref{g} imply
    \begin{equation*}
    \begin{split}
    T\left(r,\frac{1}{f-g}\right) & \leq T(r,f) + T(r,g) + O(1)\\ & \leq \frac{2}{3+\varepsilon} \sum_{k=1}^5
    \widetilde{N}_c\left(r,\frac{1}{f-a_k}\right)\\ & \leq \frac{2}{3+\varepsilon} N\left(r,\frac{1}{f-g}\right) \\
    &\leq \frac{2}{3+\varepsilon} T\left(r,\frac{1}{f-g}\right).
    \end{split}
    \end{equation*}
This is only possible when $f-g$ is a constant, say $g(z)=f(z)+k$.
But now, since $f(z)$ and $f(z)+k$ share five points out of which
at most two can be either exceptionally paired or Picard
exceptional, $k=0$, and the assertion follows.
\end{proof}

The elliptic functions $\textrm{sn}\, z$ and $1/\textrm{sn}\, z$
show that the number five cannot be replaced by four in
Theorem~\ref{5}. Namely, for both functions zero and infinity are
exceptional paired values, and they share the points $1$ and $-1$,
counting multiplicities. Therefore, $\textrm{sn}\, z$ and
$1/\textrm{sn}\, z$ share the points $-1$, $0$, $1$ and $\infty$,
ignoring pairs.

\subsection{An application to difference equations}

In this section we give an example of how to apply the obtained
results to study meromorphic solutions of difference equations. We
consider the equation
    \begin{equation}\label{diffeq2}
    w(z+1)+w(z-1)=\frac{a_2 w(z)^2+a_0}{1-w(z)^2}
    \end{equation}
where the right side is irreducible in $w$ and the coefficients
$a_j$ are constants. Equation \eqref{diffeq2} is a subcase of a
more general equation studied in \cite{halburdk:04} where it was
shown that the existence of one finite-order meromorphic solution
is sufficient to reduce a large class of difference equations into
a difference Painlev\'e equation or into a linear difference
equation, provided that the solution does not simultaneously
satisfy a difference Riccati equation. Suppose that
\eqref{diffeq2} has a finite-order meromorphic solution $w(z)$ and
consider a Laurent series expansion of $w$ in a neighborhood of a
point $z_0$ such that $w(z_0)=\delta$ with the multiplicity
$k\geq1$, where $\delta:=\pm1$. Then $w$ has a pole of order at
least $k$ at $z_0-1$ or $z_0+1$.

Consider first the case where $w(z_0+1)=\infty$ with the
multiplicity $k$ and $w(z_0-1)$ is either finite or a pole with
multiplicity strictly less than $k$. Then by iterating
\eqref{diffeq2}, we have
    \begin{equation}\label{case1}
    \begin{split}
    w(z+4n) &= \delta + \alpha(z-z_0)^k +
    O\left((z-z_0)^{k+1}\right)\\
    w(z+2n+1) &= \frac{\left((-1)^n(\frac{1}{4}n+\frac{1}{8})-\frac{1}{8}\right)(a_0+a_2)}{\alpha\delta}(z-z_0)^{-k}
    +
    O\left((z-z_0)^{1-k}\right)\\
    w(z+4n+2) &= -a_2-\delta +
    O\left((z-z_0)\right)
    \end{split}
    \end{equation}
for all $n\in\N\cup\{0\}$ and for all $z$ in a suitably small
neighborhood of $z_0$, provided that $a_2\not=0$. Since we assumed
the right side of \eqref{diffeq2} to be irreducible
$a_0+a_2\not=0$ and so $w(z+2n+1)=\infty$ for all
$n\in\N\cup\{0\}$. The iteration in the case where $w(z_0+1)$ is
finite, or a pole with low order, and $w(z_0-1)=\infty$ is
symmetric with \eqref{case1}.

Suppose now that $w(z_0)=\delta$ and $w(z_0\pm1)=\infty$ all with
the same multiplicity $k$. Then, assuming $c_1\in\C$ and
$c_{-1}\in\C$ such that $c_1c_{-1}\not=0$, we have
    \begin{equation}\label{case2}
    \begin{split}
    w(z+4n) &= \delta + \alpha(z-z_0)^k +
    O\left((z-z_0)^{k+1}\right)\\
    w(z+4n+2) &= -a_2-\delta +
    O\left((z-z_0)\right)\\
    w(z+2n+1) &= c_{2n+1}(z-z_0)^{-k}
    + O\left((z-z_0)^{1-k}\right)\\
    \end{split}
    \end{equation}
for all $n\in\Z$ as long as none of the constants $c_{2n+1}$
vanish. But if $c_{2n_0+1}=0$ for some $n_0\in\Z$ then we are back
in the situation \eqref{case1} with the starting point
$z_0+2n_0+1$ instead of $z_0-1$. Note also that a closer
inspection of the iteration in \eqref{case2} shows that
    \begin{equation}\label{cs}
    c_{k\pm4}=c_k + \frac{a_2+a_0}{2\alpha}
    \end{equation}
for all $k\in\Z$.

The final possibility is that $w(z_0)=\delta$ with the
multiplicity $k$ and $w(z_0\pm1)=\infty$ for both choices of the
sign with the multiplicity strictly greater than $k$. But in this
case it is immediately seen that $w(z)$ has a pole with the same
order in $z_0+2n+1$ for all $n\in\Z$.

We conclude that all poles, $1$-points and $-1$-points of $w$
appear in lines where each point is separated from its neighbors
by the constant $4$, with the possible exception of the endpoints
of sequences of points appearing as a part of \eqref{case1}. In
fact for our purposes it is sufficient to know that all poles and
$\delta$-points of $w$ appear in groups of four or more, with
$4$-separation. Assume that $w$ is not periodic with period four.
Then by Theorem~\ref{2nd2},
    \begin{equation*}
    \begin{split}
    T(r,w) &\leq \widetilde{N}_4(r,\infty) + \widetilde{N}_4(r,1) +\widetilde{N}_4(r,-1) +
    S(r,w)\\ &\leq
    \frac{1}{4}N(r,\infty) + \frac{1}{4}N(r,1) +\frac{1}{4}N(r,-1) +
    S(r,w)\\
    &\leq \frac{3}{4}T(r,w) + S(r,w),
    \end{split}
    \end{equation*}
which is a contradiction. Therefore, either $a_2=0$, or $w$ is
periodic with period $4$ or of infinite order.

Suppose finally that $w$ is periodic with period $4$. Then $1$ and
$-1$ are Picard exceptional values of $w$ by \eqref{case1},
\eqref{case2} and \eqref{cs}. Therefore all poles of $w$ appear in
lines where each pole is separated from its neighbors by the
constant $2$, and so $w$ is periodic with period $2$. But then, by
periodicity, $w(z+1)$, $w(z-1)$ and $w(z+1)+w(z-1)$ are infinite
simultaneously. On the other hand, the right side of
\eqref{diffeq2} is never infinite since the values $\pm1$ are
Picard exceptional. Hence also the value infinity is Picard
exceptional for $w$, and therefore $w$ is a constant by Picard's
theorem. We conclude that if \eqref{diffeq2} has a non-constant
meromorphic solution of finite order then $a_2=0$.

The existence of finite-order meromorphic solutions of
\eqref{diffeq2} is guaranteed in the case $a_2=0$, $a_0\not=0$.
Then \eqref{diffeq2} has solutions of the form
    \begin{equation}\label{sol}
    w(z)=\frac{\alpha\, \textrm{sn}(\Omega z+C)+\beta}{\gamma\, \textrm{sn}(\Omega z+C)+\delta}
    \end{equation}
where $C\in\C$ is arbitrary, and
$\alpha,\beta,\gamma,\delta,\Omega$ are certain constants
depending on another free parameter. The meromorphic solutions
\eqref{sol} are of order $2$ and periodic, but not of period $4$.

\section{Discussion}

Nevanlinna's second main theorem implies that a non-constant
meromorphic function cannot have too many points with high
multiplicity. In this study a difference analogue of the second
main theorem of Nevanlinna theory was given, which shows that a
non-periodic finite-order meromorphic function cannot have many
values which only appear in pairs, separated by a fixed constant.
Then a number of results on the value distribution of finite-order
meromorphic functions were derived by combining existing proof
techniques from Nevanlinna theory together with the difference
analogue of the second main theorem. These include analogues of
Picard's theorem, the theorem on the deficiency sum and the
theorem on meromorphic functions sharing five values. Sharpness of
these results was discussed with the help of examples. Also, an
example of how to apply some of these results to study complex
difference equations was given.

All concepts of Nevanlinna theory related to ramification have a
natural difference analogue. For instance, constant functions are
analogous to periodic functions, and a pole with multiplicity
$n>1$ is analogous to a line of $n$ poles with the same
multiplicity, each separated from its neighbors by a fixed
constant. Similarly as a pole is counted only once in the counting
function $\overline{N}(r,f)$ regardless of its multiplicity, only
one pole from the above line of poles contributes to
$\widetilde{N}_c(r,f)$. However, some notions in the difference
Nevanlinna theory seem to go, in a sense, further than their
classical counterparts. If a line of points consists of poles with
different multiplicities, the contribution from these poles to
$\widetilde{N}_c(r,f)$ is nevertheless strictly less than the
contribution to $N(r,f)$. Therefore this situation is still
exceptional in the sense of the difference deficiency relation
\eqref{defsumpair}. On the other hand, if all poles in the line
have similar enough Laurent series expansions, then the
contribution to $\widetilde{N}_c(r,f)$ from these poles may be
\textit{negative}. This implies that the maximal value two in the
difference deficiency relation \eqref{defsumpair} may be attained
by one value~$a$, which is impossible for the classical
deficiencies \eqref{defsum}.

\section{Open problems}

In addition to his ground-breaking results in the field of value
distribution theory, Nevanlinna proposed a number of problems many
of which have remained open until recently. In this section we
briefly discuss two of them.

\subsection{Inverse problem}

The \textit{inverse problem} for the deficiency relation is to
find a meromorphic function $f$ which at prescribed points has
certain non-zero deficiencies and ramification indices. This
problem was proposed and partially solved by Nevanlinna himself,
see \cite{nevanlinna:74}, but the complete solution had to wait
until 1977 when Drasin \cite{drasin:77} settled the issue by a
clever use of quasi-conformal mappings. Later on Drasin
\cite{drasin:87} established a related corollary by F. Nevanlinna,
which states that if a meromorphic function $f$ has finite order
$\lambda$ and $\sum_a\delta(a,f)=2$ then $2\lambda$ is a natural
number greater or equal to two. In the view of Corollary \ref{rel}
it is natural to ask under what conditions it is possible to find
a meromorphic function of finite order for which the pair index
$\pi(a,f)$ and the deficiency $\delta(a,f)$ have certain non-zero
values at prescribed points~$a$?

\subsection{Slowly moving targets}

Another question proposed by Nevanlinna is whether or not the
relation \eqref{defsum} remains valid if the sum is taken over all
small functions with respect to $f$. Partial answer was given by
Steinmetz \cite{steinmetz:86} and Osgood \cite{osgood:85} who
showed that
    \begin{equation*}
    \sum_a \delta(a,f) \leq 2
    \end{equation*}
where the sum is taken over distinct small functions with respect
to $f$. A complete solution to this problem was given recently by
Yamanoi \cite{yamanoi:04} who showed that \eqref{defsum} indeed
remains valid if the sum is taken over the larger field small
functions, rather than just constants. Similarly we propose that
the \eqref{relsum} remains valid even if the sum is taken over the
field $\mathcal{S}(f)$. It can be immediately seen, by a
modification of the reasoning in \cite[p. 47]{hayman:64}, that the
assertion holds for at most three functions.

\bigskip

\noindent R. G. Halburd

\smallskip

\noindent{\sc Department of Mathematical Sciences, Loughborough
University, Loughborough, Leicestershire, LE11 3TU, UK.}

\smallskip

\noindent\emph{E-mail address:} {\tt r.g.halburd@lboro.ac.uk}

\bigskip

\noindent R. J. Korhonen

\smallskip

\noindent{\sc University of Joensuu, Department of Mathematics,
P.~O.~Box 111, FIN-80101 Joensuu, FINLAND.}

\smallskip

\noindent\emph{E-mail address:} {\tt risto.korhonen@joensuu.fi}

\end{document}